# Riemann's Zeta Function.
# The alternating Xi-Function $\xi_a(s)$.

Renaat Van Malderen

July 2016

*Abstract*

As well known, the study of Riemann's zeta function ζ(s) involves the related entire function $\xi(s)$. A close relative of ζ(s) is the alternating zeta function $\eta(s)$. Similar to ζ(s), also $\eta(s)$ has a corresponding entire function $\xi_a(s)$. After establishing its definition and a related functional equation, formulas based on incomplete gamma functions are worked out, allowing to numerically evaluate $\xi_a(s)$. Examples verifying the obtained formulas are included.

1. Introduction.

The functional equation for the Riemann zeta function ζ(s) defines a related function of the complex variable s:

$$\xi(s) = \frac{s(s-1)}{2} \pi^{-\frac{s}{2}} \Gamma\left(\frac{s}{2}\right) \zeta(s) \qquad (1)$$

$\xi(s)$ is an even entire function with respect to ½, i.e.

$$\xi\left(s - \frac{1}{2}\right) = \xi\left(\frac{1}{2} - s\right) \qquad (2)$$

The present paper will mostly deal with the alternating zeta function:

$$\eta(s) = \sum_{n=1}^{\infty} (-1)^{n+1} n^{-s} \qquad (3)$$

(3) is convergent for Re(s)>0 and absolutely convergent for Re(s)>1. $\eta(s)$ may be analytically continued to the entire s plane by simply using its relation to ζ(s):

$$\eta(s) = (1 - 2^{1-s}) \zeta(s) \qquad (4)$$

Unlike ζ(s) which has a simple pole at s=1, $\eta(s)$ is an entire function. As easily verified by (4), all the zeros of ζ(s) (both trivial and non trivial) are also zeros of $\eta(s)$. Due to the factor $(1-2^{1-s})$, $\eta(s)$ has additional simple zeros for

$$s = 1 \pm \left(\frac{2n\pi i}{\ln 2}\right) \text{ with n=1,2,3,...} \qquad (5)$$

For s=1, $\eta(s)$ remains finite and from (3):

$$\eta(1) = \ln 2 \qquad (6)$$

while $\eta(0) = (1 - 2^{1-0})\zeta(0) = (-1)\left(-\frac{1}{2}\right) = \frac{1}{2}$ (7)

As shown below, it is possible to define for $\eta(s)$ an entire even function which we will call the Alternating Xi function $\xi_a(s)$ and which plays a similar role for $\eta(s)$ as $\xi(s)$ for $\zeta(s)$.

## 2. The alternating Xi Function $\xi_a(s)$.

Consider the integral $I_n = \int_0^\infty E(-\pi n^2 x) x^{\frac{s}{2}-1} dx$ with Re(s)>0.
Through substitution $I_n$ may be expressed as:

$$I_n = \pi^{-\frac{s}{2}} n^{-s} \Gamma\left(\frac{s}{2}\right) \qquad (8)$$

We define: $\phi(x) = \sum_{n=1}^\infty (-1)^{n+1} E(-\pi n^2 x)$ (9)

The series (9) is absolutely and uniformly convergent for x>0.

(9) $\Rightarrow \pi^{-s/2} \Gamma\left(\frac{s}{2}\right) \eta(s) = \int_0^\infty \phi(x) x^{\frac{s}{2}-1} dx$ (10)

Multiplying (10) by $2^s$:

$$2^s \pi^{-\frac{s}{2}} \Gamma\left(\frac{s}{2}\right) \eta(s) = \int_0^\infty \phi(x)(4x)^{\frac{s}{2}-1} dx = \int_0^\infty \phi\left(\frac{x}{4}\right)(x)^{\frac{s}{2}-1} dx$$

$$(1 - 2^s) \pi^{-\frac{s}{2}} \Gamma\left(\frac{s}{2}\right) \eta(s) = \int_0^\infty \left(\phi(x) - \phi\left(\frac{x}{4}\right)\right)(x)^{\frac{s}{2}-1} dx \qquad (11)$$

We define: $\varphi(x) = \phi(x) - \phi\left(\frac{x}{4}\right)$ (12)

(12) is of course also absolutely and uniformly convergent for x>0. As demonstrated in addendum -a, the following property for $\varphi(x)$ holds:

$\varphi\left(\frac{1}{x}\right) = \varphi(x)\sqrt{x}$ (13)

Using (12), $\varphi(0)$ is strictly undefined. But $\lim_{x\to\infty} \varphi(x) = 0$. Using (13) also $\lim_{x\to 0} \varphi(x) = 0$.
As explained later in this paper (see formula (17)), $\varphi(x)$ may be written as
$\varphi(x) = \sum_{m=0}^\infty T_m(x)$ with:

$$T_m(x) = -E\left(\frac{-\pi x}{4}(4m+1)^2\right) + 2E\left(\frac{-\pi x}{4}(4m+2)^2\right) - E\left(\frac{-\pi x}{4}(4m+3)^2\right)$$

For large enough x-values the first term of $T_m(x)$ is dominant and renders $T_m(x)$ negative. (Already for x=1 and m=0, we have $T_0(1)$= -0.37036….). For larger x values and larger m, this remains true and more pronounced. This implies that for all m and at least x≥1, all $T_m(x)<0$ and therefore $\varphi(x) < 0$. (13) then implies that $\varphi(x) < 0$ for all x>0. Differentiating (13) and putting x=1 yields:

$$\varphi'(1) = -\varphi(1)/4 \qquad (14)$$

(14) $\Rightarrow \varphi'(1) > 0$ Numerical evaluation of (12) yields:
$\varphi(1) = -0.370361\ldots$ $\qquad\qquad\qquad \varphi'(1) = +0.092590\ldots$
$|\varphi(x)|max = |\varphi(0.8666\ldots)|=0.377066$

To calculate $\varphi(x)$ for x values close to zero, it is advantageous to use (13), instead of (12) the latter requiring to extend the calculation over several m values. E.g. $\varphi(0.1) \cong \varphi(10)\sqrt{10}$ = -0.00122…. and $\varphi(0.05) \cong \varphi(20)\sqrt{20}$ = -6.7E(-7).

We define the alternating Xi function by:

$$\xi_a(s) = (1 - 2^s)\pi^{-s/2}\Gamma\left(\frac{s}{2}\right)\eta(s) = \int_0^\infty \varphi(x)x^{\frac{s}{2}-1}dx \qquad (15)$$

Please notice $\xi_a(s)$ is the Mellin transform $M(s/2)$ of $\varphi(x)$. The integral (15) is convergent for Re(s)>0. As we will indicate further on, its validity stretches over the entire s plane. Also from (15) it may be seen that $\xi_a(s)$ is an entire function: The pole of $\Gamma\left(\frac{s}{2}\right)$ at s=0 is cancelled by the corresponding zero of $(1 - 2^s)$. The other poles of $\Gamma\left(\frac{s}{2}\right)$ at s=-2,-4,-6,… are cancelled by the zeros of $\eta(s)$ for the same s –values.

Table-1 gives values for $(1 - 2^s), \zeta(s), \eta(s)$, and $\xi_a(s)$ for some key s values.

Table -1

| s | $(1 - 2^s)$ | $\zeta(s)$ | $\eta(s)$ | $\xi_a(s)$ |
|---|---|---|---|---|
| 0 | 0 | -1/2 | 1/2 | -ln(2) |
| 1 | -1 | $+\infty$** | ln(2) | -ln(2) |
| $\pm 2n\pi i/ln2$* | 0 | $\neq 0$ | $\neq 0$ | 0 |
| $1 \pm 2n\pi i/ln2$* | -1 | $\neq 0$ | 0 | 0 |

*n=1,2,3,… $\qquad\qquad$ ** for s=1+

We now demonstrate that $\xi_a(s)$ is an even function with respect to s=1/2. We split $\int_0^\infty$ in (15) into $\int_0^1 + \int_1^\infty$ and apply (13) to $\int_0^1$ :

$$\xi_a(s) = \int_0^1 \frac{\varphi\left(\frac{1}{x}\right)}{\sqrt{x}} x^{\frac{s}{2}-1}dx + \int_1^\infty \varphi(x)x^{\frac{s}{2}-1}dx$$

As an intermediate step we substitute y for 1/x and then returning to x again we have:

$$\xi_a(s) = \int_1^\infty \varphi(x)\frac{dx}{x}\left[x^{\frac{s}{2}} + x^{\frac{1-s}{2}}\right] \qquad (16)$$

(16) $\Rightarrow \xi_a(s) = \xi_a(1 - s)$. This relation together with (15) yields a <u>functional equation</u> for $\eta(s)$:

$$\eta(s) = \frac{2^s - 2}{1 - 2^s} \pi^{s-1} \Gamma(1-s) \sin\left(\frac{\pi s}{2}\right) \eta(1-s) \quad (16-a)$$

Apart from "trivial" zeros introduced by both factors $(1 - 2^s)$ and $\eta(s)$ in (15), $\xi_a(s)$ has the same "nontrivial" zeros as $\zeta(s)$.

Value of $\xi_a\left(s = \frac{1}{2}\right)$: In reference [5, p13, using formula (28)] $\xi\left(\frac{1}{2}\right)$ was calculated to be $\xi\left(\frac{1}{2}\right) = 0.4971208\ldots$
Inserting this result into formula (28) of this paper yields:
$\xi_a\left(\frac{1}{2}\right) = -0.6823392\ldots$ and via (15):
$\eta\left(\frac{1}{2}\right) = 0.6048986\ldots$ and $\zeta\left(\frac{1}{2}\right) = -1.4603544\ldots$

Some remarkable cases of formula (15):
Since $\xi_a(s) = \xi_a(1-s)$ we have : $\int_0^\infty \varphi(x) x^{\frac{s}{2}-1} dx = \int_0^\infty \varphi(x) x^{-\frac{s}{2}-\frac{1}{2}} dx$.
This yields for s=0: $\int_0^\infty \varphi(x) \frac{dx}{x} = \int_0^\infty \varphi(x) \frac{dx}{\sqrt{x}} = -ln2$ (see table -1).
$\xi_a(2)$ may be obtained from (15); Noting that $\Gamma\left(\frac{2}{2}\right) = 1$ and $\eta(2) = \frac{\pi^2}{12}$, we have
$\int_0^\infty \varphi(x) dx = -\frac{\pi}{4}$.

The wording of H.M. Edwards used for the case of $\xi(s)$ [2, p.16] is also applicable to (16): Because $\varphi(x)$ decreases more rapidly than any power of x as x→∞, the integral in (16) converges in the entire s plane.

### 3. Evaluation of $\xi_a(s)$ via Incomplete Gamma Functions.

The function $\varphi(x)$ given in (12) equals the difference of two infinite series of which the terms are of the type:

$(-1)^{n+1} E(-\pi n^2 x)$ and $(-1)^{n+1} E(-\pi n^2 \frac{x}{4})$. As easily verified, the total expression may be written as:

$$\varphi(x) = \sum_{m=0}^\infty \left[-E\left(\frac{-\pi x}{4}(4m+1)^2\right) + 2E\left(\frac{-\pi x}{4}(4m+2)^2\right) - E\left(\frac{-\pi x}{4}(4m+3)^2\right)\right] \quad (17)$$

Let: $A_m = \frac{\pi}{4}(4m+1)^2$, $B_m = \frac{\pi}{4}(4m+2)^2$, $C_m = \frac{\pi}{4}(4m+3)^2$

$$\varphi(x) = \sum_{m=0}^\infty [-E(-A_m x) + 2E(-B_m x) - E(-C_m x)] \quad (18)$$

The first integral in (16) then becomes:

$$\int_1^\infty \frac{dx}{x} x^{\frac{s}{2}} \sum_{m=0}^\infty [-E(-A_m x) + 2E(-B_m x) - E(-C_m x)] \quad (19)$$

We introduce new variables $\lambda_A = A_m x$, $\lambda_B = B_m x$, $\lambda_C = C_m x$, with corresponding integration limits: $A_m \to \infty$; $B_m \to \infty$; $C_m \to \infty$.

Let's consider for a given m in (19) the resulting integral for the term in $A_m$:

$$\int_1^\infty dx\, x^{\frac{s}{2}-1} E(-A_m x) = \frac{\int_{A_m}^\infty d\lambda_A\, \lambda_A^{\frac{s}{2}-1} E(-\lambda_A)}{A_m^{s/2}} = \frac{\Gamma(\frac{s}{2}, A_m)}{A_m^{s/2}}$$

in which $\Gamma\left(\frac{s}{2}, A_m\right)$ stands for the upper incomplete gamma function [3, p.260], [4].

Proceeding similarly for the other terms in (16) and (18):

$$\xi_a(s) = \sum_{m=0}^\infty \left[ -\frac{\Gamma\left(\frac{s}{2}, A_m\right)}{A_m^{\frac{s}{2}}} + \frac{2\Gamma\left(\frac{s}{2}, B_m\right)}{B_m^{\frac{s}{2}}} - \frac{\Gamma\left(\frac{s}{2}, C_m\right)}{C_m^{\frac{s}{2}}} - \frac{\Gamma\left(\frac{1-s}{2}, A_m\right)}{A_m^{\frac{(1-s)}{2}}} + \frac{2\Gamma\left(\frac{1-s}{2}, B_m\right)}{B_m^{\frac{(1-s)}{2}}} \right.$$
$$\left. - \frac{\Gamma\left(\frac{1-s}{2}, C_m\right)}{C_m^{\frac{(1-s)}{2}}} \right] \quad (20)$$

In particular for the critical line (s=1/2+it) and in order to keep formulas short we put:

$$\omega = \frac{\frac{1}{2}+it}{2} = \frac{1}{4} + \frac{it}{2} \qquad \bar{\omega} = \frac{1-(\frac{1}{2}+it)}{2} = \frac{1}{4} - \frac{it}{2} \quad (21)$$

$$\xi_a\left(\frac{1}{2}+it\right) = \sum_{m=0}^\infty \left[ -\frac{\Gamma(\omega, A_m)}{A_m^\omega} + \frac{2\Gamma(\omega, B_m)}{B_m^\omega} - \frac{\Gamma(\omega, C_m)}{C_m^\omega} - \frac{\Gamma(\bar\omega, A_m)}{A_m^{\bar\omega}} + \frac{2\Gamma(\bar\omega, B_m)}{B_m^{\bar\omega}} \right.$$
$$\left. - \frac{\Gamma(\bar\omega, C_m)}{C_m^{\bar\omega}} \right] \quad (22)$$

Or:

$$\xi_a\left(\frac{1}{2}+it\right) = 2\mathrm{Re} \sum_{m=0}^\infty \left[ -\frac{\Gamma(\omega, A_m)}{A_m^\omega} + \frac{2\Gamma(\omega, B_m)}{B_m^\omega} - \frac{\Gamma(\omega, C_m)}{C_m^\omega} \right] \quad (23)$$

Alternative expression for $\xi_a(s)$:

By splitting (15) again into two integrals and this time converting $\int_1^\infty$ into $\int_0^1$ we obtain:

$$\xi_a(s) = \int_0^1 \varphi(x) \frac{dx}{x}\left[x^{\frac{s}{2}} + x^{\frac{1-s}{2}}\right] \quad (24)$$

By using (18) into (15) we have:

$$\xi_a(s) = \int_0^\infty \varphi(x) x^{\frac{s}{2}-1} dx = \Gamma\left(\frac{s}{2}\right) \sum_{m=0}^\infty \left[-A_m^{-s/2} + 2B_m^{-s/2} - C_m^{-s/2}\right] \quad (25)$$

Eliminating $\Gamma\left(\frac{s}{2}\right)$ out of the above as well as in (15):

$$\eta(s)(1 - 2^{-s}) = \sum_{m=0}^\infty [(4m+1)^{-s} - 2(4m+2)^{-s} + (4m+3)^{-s}] \quad (25-a)$$

It is easily checked that (25-a) equates to $\eta(s) = \sum_{n=1}^\infty (-1)^{n+1} n^{-s}$ for Re(s)>0.

By using the relation $\Gamma(a) = \gamma(a,z) + \Gamma(a,z)$ between the lower and upper incomplete gamma functions in (20) and taking (25) into account and noting that $\xi_a(s) = \xi_a(1-s)$ we end up with the somewhat surprising result:

$$\xi_a(s) = \sum_{m=0}^\infty \left[ -\frac{\gamma\left(\frac{s}{2}, A_m\right)}{A_m^{\frac{s}{2}}} + \frac{2\gamma\left(\frac{s}{2}, B_m\right)}{B_m^{\frac{s}{2}}} - \frac{\gamma\left(\frac{s}{2}, C_m\right)}{C_m^{\frac{s}{2}}} - \frac{\gamma\left(\frac{1-s}{2}, A_m\right)}{A_m^{\frac{(1-s)}{2}}} + \frac{2\gamma\left(\frac{1-s}{2}, B_m\right)}{B_m^{\frac{(1-s)}{2}}} \right.$$
$$\left. - \frac{\gamma\left(\frac{1-s}{2}, C_m\right)}{C_m^{\frac{(1-s)}{2}}} \right] \quad (26)$$

## 4. A Numerical Example.

In [5, p. 13] the value of the Xi function $\xi(s)$ was computed for a point $s = \frac{1}{2} + 12i$ on the critical line: $\xi(0.5 + 12i) = 0.008823639 \ldots$ (27)

In the present example we calculate $\xi_a(s)$ for the same s value using (23). We will then compare the result with the value for $\xi(0.5 + 12i)$ as given in (27) after converting $\xi_a(0.5 + 12i)$ to the corresponding $\xi(s)$ value using the formula:

$$\xi(s) = \frac{s(s-1)\xi_a(s)}{2(1-2^s)(1-2^{1-s})} \quad (28)$$

In our example the $\omega$-value used in (23) equals:

$$\omega = \frac{1}{4} + \frac{it}{2} = \frac{1}{4} + 6i$$

As will become clear, the summation $\sum_{m=0}^\infty$ for this case may be limited to m=0. Already for m=1 the relative value of the numbers involved is negligible. Table -2 shows the numbers. Γ(ω) was obtained using Stirling's formula. The values for Γ(ω,a) were obtained by first computing γ(ω,a) via series expansion [5, p. 12, formula 38,39] and then taking the difference:

Γ(ω,a) = Γ(ω) - γ(ω,a)

As obvious from table -2, Γ(ω,a) values for $a = \frac{25\pi}{4}, 9\pi, \frac{49\pi}{4}$ may be neglected. Also unsurprisingly, for these larger a-values the γ(ω,a) approach increasingly Γ(ω) itself, resulting in very small differences with Γ(ω,a) .

| Table -2 | $\omega = \frac{1}{4} + 6i$ | Γ(ω,a) = Γ(ω) - γ(ω,a) |
|---|---|---|
| m | | Γ(ω)= -0.000044667603-0.000121313951i |
| 0 | | γ(ω,π/4)= -0.072862357-0.002369978i |
| | | γ(ω,π)=+0.000522089-0.0091831347i |
| | | γ(ω,9π/4)= -0.000191090527-0.000092456237i |
| 1 | | γ(ω,25π/4)= -0.000044667841-0.000121313756i |
| | | γ(ω,9π)= -0.000044667616-0.000121313949i |
| | | γ(ω,49π/4)= -0.000044667616-0.000121313949i |
| 0 | | Γ(ω,π/4)= +0.072817689397+0.00224866405i |
| | | Γ(ω,π)=-0.000566756603+0.00906182075i |
| | | Γ(ω,9π/4)= +0.000146422924-0.000028857713i |
| 1 | | Γ(ω,25π/4)= +0.000000000238-0.000000000194i |
| | | Γ(ω,9π)= +0.000000000013-0.000000000001i |
| | | Γ(ω,49π/4)= +0.000000000013-0.000000000001i |

Limiting our calculations to m=0 and abbreviating:

$$\alpha_0 = 2\text{Re}\left[\frac{\Gamma\left(\omega,\frac{\pi}{4}\right)}{\left(\frac{\pi}{4}\right)^\omega}\right] = 0.013993985486\ldots$$

$$\beta_0 = 2\text{Re}\left[\frac{\Gamma(\omega,\pi)}{(\pi)^\omega}\right] = 0.00680962358\ldots$$

$$\gamma_0 = 2\text{Re}\left[\frac{\Gamma\left(\omega,\frac{9\pi}{4}\right)}{\left(\frac{9\pi}{4}\right)^\omega}\right] = 0.000147065423\ldots$$

$\xi_a(0.5 + 12i) = -\alpha_0 + 2\beta_0 - \gamma_0 = -0.000521803749\ldots$

Using (27) and (28): $\xi(0.5 + 12i) = +0.008823639811\ldots$
<u>in line with the value obtained in</u> [5, p. 13].

Using (15) with $s = \frac{1}{2} + 12i$, we may obtain the value for $\eta = \frac{1}{2} + 12i$ as well as for the corresponding $\zeta(s) = \frac{\eta(s)}{(1-2^{1-s})}$.

(15) $\Rightarrow \eta(s) = \frac{\xi_a(s)\pi^{s/2}}{(1-2^s)\Gamma\left(\frac{s}{2}\right)}$      (29)

This yields: $\eta(0.5 + 12i) = 2.601080675 + 0.0684891589i$

$\zeta(0.5 + 12i) = +1.015935342 - 0.7451116651i\ldots$

## 5. Addendum –a: Demonstration of formula (13):

$$\varphi\left(\frac{1}{x}\right) = \varphi(x)\sqrt{x} \qquad (1-a)$$

with $\varphi(x) = \phi(x) - \phi\left(\frac{x}{4}\right)$ $\qquad \varphi\left(\frac{1}{x}\right) = \phi\left(\frac{1}{x}\right) - \phi\left(\frac{1}{4x}\right)$ $\qquad (2-a)$

where $\phi(x) = \sum_{n=1}^{\infty}(-1)^{n+1}E(-\pi n^2 x)$ $\qquad \phi\left(\frac{x}{4}\right) = \sum_{n=1}^{\infty}(-1)^{n+1}E\left(-\pi n^2 \frac{x}{4}\right)$ $\qquad (3-a)$

In [1, p. 124, Exercise 17] the relation $\psi\left(\frac{1}{x}\right) = \sqrt{x}\,\psi(x)$ $\quad (4-a)$ for the function $\psi(x) = \sum_{n=-\infty}^{\infty} E(-\pi n^2 x)$ may be obtained by considering the complex integral:

$$\oint \frac{E(-\pi z^2 x) dz}{E(2\pi i z) - 1} \qquad (5-a)$$

around a rectangle whose corners are $\pm(N+\tfrac{1}{2})\pm i$ with the integer N→∞. Riemann [2, p.15] in his second proof of the Zeta functional equation used a derived form of (4-a), i.e.

$$\frac{1 + 2\psi_*(x)}{1 + 2\psi_*\left(\frac{1}{x}\right)} = \frac{1}{\sqrt{x}} \qquad (6-a)$$

where $\psi_*(x)$ has a slightly different meaning i.e. $\psi_*(x) = \sum_{n=1}^{\infty} E(-\pi n^2 x)$. (6-a) is a form of the Theta functional equation studied by Jacobi. Similar to (5-a), relation (1-a) may be derived using the same complex contour as in (5-a) after however modifying (5-a) into:

$$\oint \frac{E(-\pi z^2 x) E(i\pi z) dz}{E(2\pi i z) - 1} \qquad (7-a)$$

allowing thereby for the alternating signs in (3-a). Further straightforward elaboration of (7-a) eventually yields (1-a). A more direct way to obtain (1-a) uses the transformation formulas for Thetaseries as in e.g. [6, pp. 345,346]. In this reference three Theta series are considered:

$$\left.\begin{array}{l}\vartheta(z) = \displaystyle\sum_{n=-\infty}^{\infty} E(i\pi n^2 z) \\[1em] \tilde{\vartheta}(z) = \displaystyle\sum_{n=-\infty}^{\infty} (-1)^n E(i\pi n^2 z) \\[1em] \widetilde{\widetilde{\vartheta}}(z) = \displaystyle\sum_{n=-\infty}^{\infty} E\left(i\pi\left(n+\frac{1}{2}\right)^2 z\right)\end{array}\right\} \qquad (8-a)$$

By putting z=ix, we find respectively:

$$\phi(x) = \frac{1 - \tilde{\vartheta}(ix)}{2}$$

with $\tilde{\vartheta}(ix) = \sum_{n=-\infty}^{\infty}(-1)^n E(-\pi n^2 x)$     $(9-a)$

$$\phi\left(\frac{x}{4}\right) = \frac{\tilde{\tilde{\vartheta}}(ix)}{2} + \frac{1}{2} - \vartheta(ix)$$

with $\begin{cases} \vartheta(ix) = \sum_{n=-\infty}^{\infty} E(-\pi n^2 x) & (10-a) \\ \\ \tilde{\tilde{\vartheta}}(ix) = 2\sum_{n=0}^{\infty} E\left(\frac{-(2n+1)^2 x}{4}\right) & (11-a) \end{cases}$

$$\varphi(x) = \phi(x) - \phi\left(\frac{x}{4}\right) = \frac{1}{2} - \frac{\tilde{\vartheta}(ix)}{2} - \frac{\tilde{\tilde{\vartheta}}(ix)}{2} - \frac{1}{2} + \frac{\vartheta(ix)}{2}$$

or

$$\varphi(x) = \frac{1}{2}\left(\vartheta(ix) - \tilde{\vartheta}(ix) - \tilde{\tilde{\vartheta}}(ix)\right) \qquad (12-a)$$

Verification of (9-a), (10-a), (11-a), (12-a) using (2-a) and (8-a) requires some algebraic but straightforward work. The Theta formulas as given in [6] are:

$$\left.\begin{array}{r}\vartheta\left(\frac{i}{x}\right) = \sqrt{x}\,\vartheta(ix) \\ \\ \tilde{\vartheta}\left(\frac{i}{x}\right) = \sqrt{x}\,\tilde{\tilde{\vartheta}}(ix) \\ \\ \tilde{\tilde{\vartheta}}\left(\frac{i}{x}\right) = \sqrt{x}\,\tilde{\vartheta}(ix)\end{array}\right\} \qquad (13-a)$$

Plugging (13-a) into (12-a) yields $\varphi\left(\frac{1}{x}\right) = \varphi(x)\sqrt{x}$

Renaat Van Malderen
Address: Maxlaan 21, B-2640 Mortsel, Belgium
The author can be reached by email: hans.van.malderen@telenet.be